\RequirePackage{fix-cm}
\documentclass{svjour3}                     
\smartqed  
\usepackage{graphicx}
\usepackage{amsmath} 
\usepackage{subfigure}
\begin{document}

\title{An alternative proof for Euler rotation theorem
}
%



\author{Toby Joseph         
}


\institute{T. Joseph \at
               BITS Pilani, K. K. Birla Goa Campus, Goa, India\\
              Tel.: +91-832-2580145\\
              \email{toby@goa.bits-pilani.ac.in}           
}

\date{Received: date / Accepted: date}

\maketitle

\begin{abstract}
Euler's rotation theorem states that any reconfiguration of a rigid body with one of its points fixed is equivalent 
to a single rotation about an axis passing through the fixed point. The theorem forms the basis for Chasles' theorem which
states that it is always possible to represent the general displacement of a rigid body by a translation and a rotation
about an axis. Though there are many ways to achieve this, the direction of the rotation axis and angle of rotation are independent 
of the translation vector. The theorem is important in the study of rigid body dynamics. There are various proofs available for these theorems, 
both geometric and algebraic. A novel geometric proof of Euler rotation theorem is presented here which makes use of two successive rotations 
about two mutually perpendicular axis to go from one configuration of the rigid body to the other with one of its points fixed. 

\end{abstract}

\maketitle

\section{Introduction}
\label{sec1}
Study of rigid body dynamics is one of the important topics in classical mechanics. As in the case of point particle
dynamics, a good understanding of the kinematic description of the rigid body motion is essential for studying its dynamics.
One of the crucial results related to the kinematic description of motion of a rigid body is Chasles' theorem\cite{mozzi,chasles,jackson}. 
It states that a general displacement of a rigid body can be described by a translation and a rotation about an
axis. Further, even though the translation vector is not unique, the orientation of axis of rotation and the 
angle of rotation will be the same for going from one configuration to the other \cite{whittaker,banach,pars}.
This natural split of general motion into translational and rotational parts allows for studying the dynamics of translation and rotation separately. 
The former reduces to a point particle like dynamics and the latter is described using Euler equations in rigid body dynamics. It turns out that
one can use the freedom in choosing the translation vector so as to make its direction coincide with that of axis of rotation. This is the
content of the Mozzi-Chasles' theorem and is central to the study of dynamics of rigid body using screw theory\cite{ball}. 
Screw theory is extensively used in the modern day study of robotics \cite{davidson,murray}.

Chasles' theorem is an extension of an earlier theorem due to Euler, referred to as 
Euler's rotation theorem. Euler's rotation theorem states that any reconfiguration of a rigid body with one of its points fixed is equivalent 
to a single rotation about an axis passing through the fixed point. In other words, whatever way a sphere might be rotated around its center, 
a diameter can always be chosen whose direction in the rotated configuration would coincide with that in the original configuration.
The original proof given by Euler himself is a geometrical one\cite{euler,palais1}. The proof looks at the initial and final configuration of a great
circle on the sphere and gives the recipe to construct a point, which is subsequently proved to be the point through which the
axis of rotation passes.

There are various other proofs available for the rotation theorem both geometric \cite{banach,pars} as well as algebraic ones \cite{palais1,palais2}. 
The algebraic proofs typically require a familiarity with rotation matrices and their properties or with ideas from group theory. The most
commonly seen analytic proof \cite{beatty,goldstein,thurnauer} uses the orthogonality property of three-dimensional rotation matrices
to show that they always have an eigenvector with eigenvalue equal to one. This proof makes use of the result that eigenvalues of an orthogonal matrix
have modulus one. Another one of the geometric proofs \cite{banach} 
involves looking at the displacement of a segment under the rotations and constructs planes of symmetry using the end points of the original and displaced segments.
The intersection of these symmetry planes is then shown to be the axis of rotation. This is by far the most transparent of the existing proofs.
The proof for Euler rotation theorem given by Pars \cite{pars} involves going from the initial configuration to the final configuration using two rotations: 
the first one is a rotation by angle $\pi$ about the center of great circle connecting one of the original points and its final location and a second
rotation about an axis passing through this final location. The invariant point is then established by a construction.

Aim of the current work is to give a geometric proof of Euler rotation theorem that is different from the existing ones. We shall first develop the 
key ideas involved in the description of rigid body and use this in setting up the proof of the theorem. After proving Euler's rotation theorem, 
which is crux of the paper, Chasles' theorem is derived. The plan of the paper is as follows: In section I, we derive the number 
of degrees of freedom of a rigid body. Further, we set up a scheme for describing the general displacement of the rigid body which we make use of
subsequently. Euler's rotation theorem is derived in section II and the proof for Chasles' theorem is given in section III. 
To make the scope of discussion wider and complete, we shall then look at a few consequences of the theorems proved involving the idea of 
screw axis and rigid body motion in two dimensions. 

\section{Description of rigid body displacement}
\label{sec2}
Rigid body is defined as a collection of particles whose mutual distances remains invariant. In three dimensions, $N$ independent particles have
$3N$ degrees of freedom. But if the particles constitute a rigid body, the degrees of freedom is reduced to $6$ (for the case when $N \ge 3$).
This is so because the constraint equations that comes from the invariant inter-particle separations makes $3N - 6$ of the original $3N$ variables
dependent on the remaining $6$. Let us prove this result rigorously. 

\begin{figure}
\begin{center}
\includegraphics[width=0.9\linewidth]{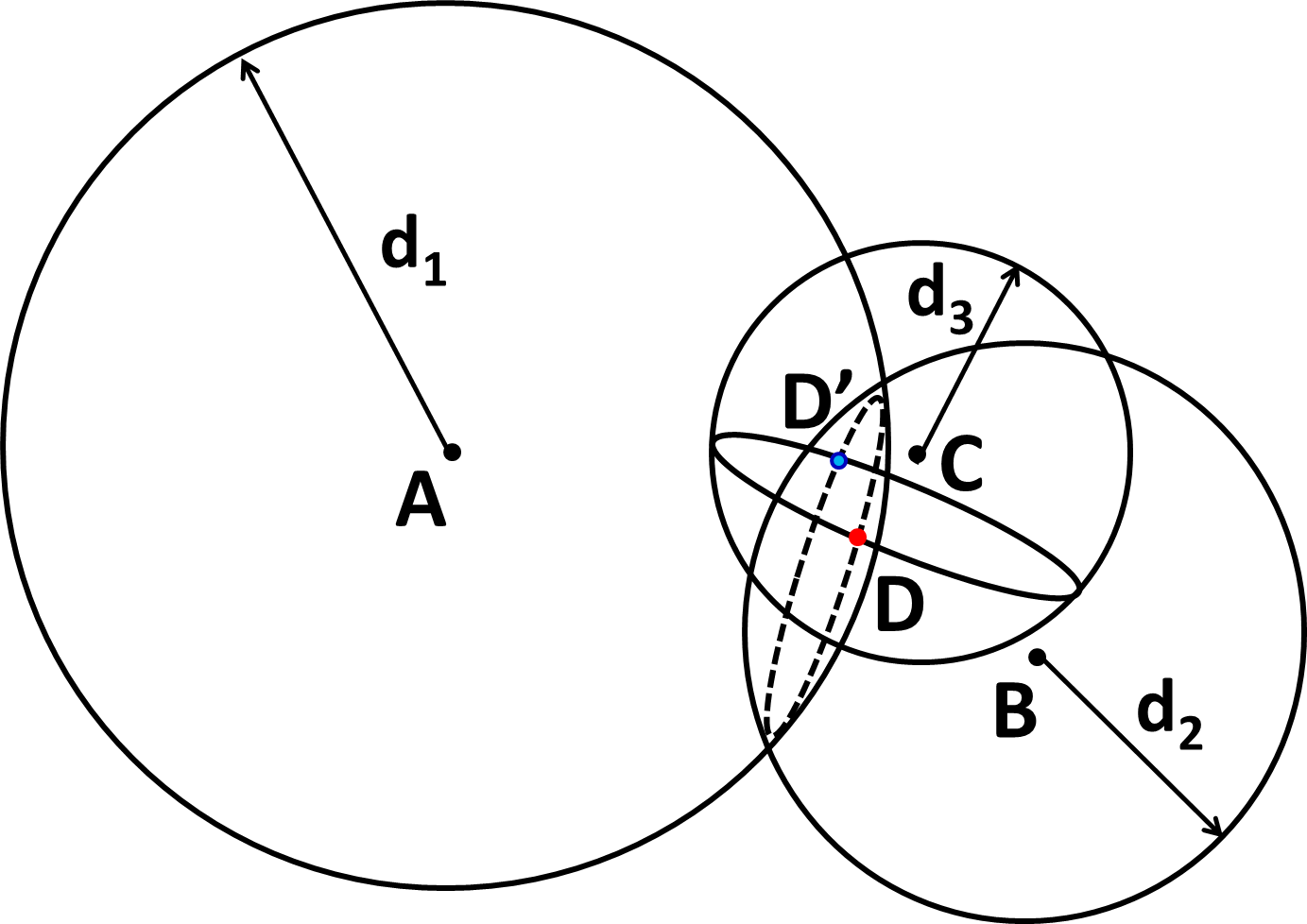} 
\caption{The construction for finding out the number of degrees of freedom of a rigid body in $3$ dimensions. If three particles $A$, $B$ and $C$ are
fixed, then there are only two possible locations for the fourth particle whose distance from the other there are fixed by rigid body constraints. The possible
locations for the fourth particle $D$ are shown as $D$ (red dot) and $D'$ (blue dot) in the figure. They are related by a reflection about the plane containing
$A$, $B$ and $C$ particles. The dashed circle is the intersection of the spheres centered around $A$ and $B$ (referred to as circle $S$ in the text). 
Points $D$ and $D'$ are the intersections of this circle with the sphere centered at $C$.}
\label{fig02}
\end{center}
\end{figure}

We first show that a rigid body configuration is completely defined once coordinates 
of any three non-collinear particles in the rigid body are specified. Assume that positions of particles (named $A, B$ and $C$) are given. Consider now
a fourth particle, $D$. Since the body is rigid, distance between particles $A$ and $D$ (say, $d_{1}$), $B$ and $D$ ($d_{2}$) and between $C$ and $D$
($d_{3}$) are fixed. Construct a sphere of radius $d_{1}$ centered around particle $A$ as shown in Fig. \ref{fig02}. It is clear that particle $D$ has to be residing on the surface
of this sphere. Now construct a second sphere of radius $d_{2}$ centered around particle $B$. The rigidity constraint will imply that particle $D$ has to
lie on the circle (call it $S$) formed by the intersection of these two spheres. Note that if the spheres do not intersect, the constraints will not be consistent with
that of a rigid body configuration. A third sphere of radius $d_{3}$ centered at $C$ will intersect the circle $S$ at two points ($D$ and $D'$ in Fig. \ref{fig02}) implying that, 
once the three points are fixed a fourth particle can only be placed at two possible points consistent with the rigidity constraints. The two possible points are related
to each other by a reflection about the plane containing particles $A, B$ and $C$ and corresponds to the mirror images of each other. So if we put an
additional constraint that the handedness of rigid body is preserved there is only a unique position where particle $D$ can be placed and hence does 
not require any further coordinates to be specified. But particle $D$ is completely arbitrary and could be any of the particles in the rigid body other than
the original three particles whose positions were specified. 

Thus we see that the number of degrees of freedom of a rigid body in three dimensions is the
same as that of a rigid body containing three non-collinear particles. Three independent particle have $9$ degrees of
freedom. Since the body is rigid, there are $3$ constraint equations specifying the mutual separation between the $3$ particles. 
The difference gives the degrees of freedom of a rigid body to be $6$. 

The next question we address is about how to describe the displacement of the rigid body from one configuration to another. There are multiple ways
to do this. We shall adopt a scheme which will make it convenient for us to construct the proof for Chasles' theorem. Consider two configurations $I$ and $II$
of the rigid body. Consider three points $P_1, P_2$ and $P_3$ (which are non-collinear) in the rigid body in configuration I. In the final configuration let these 
points move over to locations $P_1', P_2'$ and $P_3'$ respectively. To go from the configuration
$I$ to configuration $II$, we will carry out the following three steps:
\begin{enumerate}
\item  Translate the body by vector $\overrightarrow{P_1P_1'}$. This ensures that the point $P_1$ is in its final position $P_1'$.
\item Let $P_2''$ be the location of points $P_2$  after this translation. 
Consider the plane formed by the vectors $\overrightarrow{P_1' P_2''}$ and $\overrightarrow{P_1' P_2'}$. This is the equatorial plane shown in Fig. \ref{fig03}. Rotate the
rigid body about an axis perpendicular to this plane and passing through the point at $P_1'$ such that the point $P_2''$ is in its final position $P_2'$. Note that if
$P_2''$ is the same as $P_2'$ then this step need not be carried out. 
\item Let $P_3''$ be the location of the original point $P_3$ after the above two operations. The rigid body can now be rotated about the axis passing through $P_1'$ 
and $P_2'$ (see Fig. \ref{fig03}) such that the point $P_3''$ is in its final position $P_3'$.  
\end{enumerate}
These steps will ensure that the rigid body has been displaced to its final configuration.

\section{Proof of Euler rotation theorem}
\label{sec3}
\begin{figure}
\begin{center}
\includegraphics[width=0.7\linewidth]{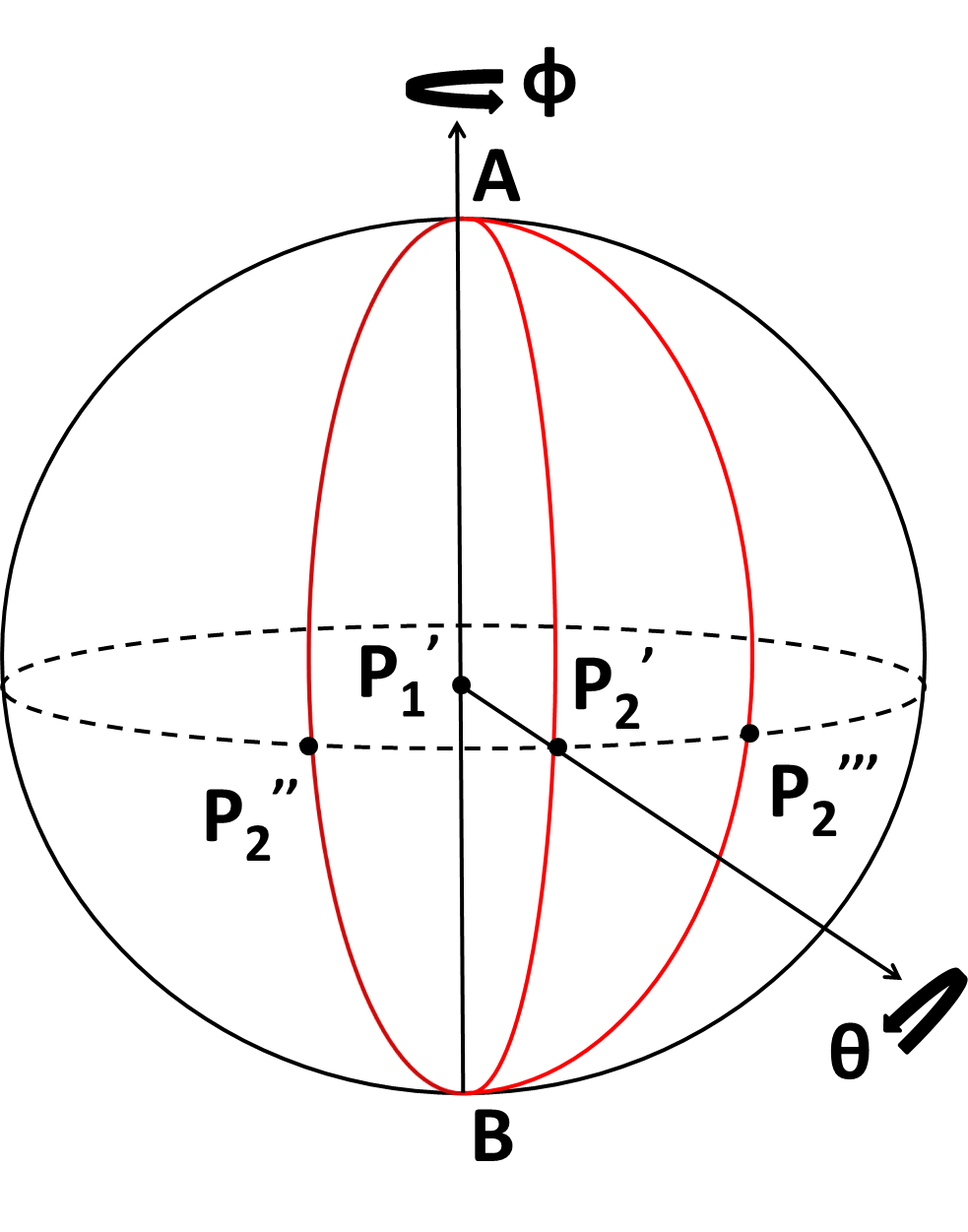} 
\caption{Steps $2$ and $3$ of the procedure for going from one configuration to the other of the rigid body. The point $P_1'$ is in its final position after step-1.
Rotation about $AB$ axis (which is perpendicular to the vectors connecting $P_1'$ to $P_2''$ and $P_1'$ to $P_2'$) by amount $\phi$ carries $P_2''$ to its final position. 
The second rotation by angle $\theta$ about $P_1' P_2'$ will bring the rigid body to its final configuration. Note that we have oriented the figure such that the $AB$ axis is
vertical. We have also not shown $P_3''$, or $P_3'$ in the figure. These points in general will not lie on the surface of the sphere shown.}
\label{fig03}
\end{center}
\end{figure}

We are now in a position to prove Euler rotation theorem. In order to prove the theorem let us look at the steps $2$ and $3$ in the scheme described above 
to go from one configuration to another of a rigid body. Note that after step $1$ the point $P_1$ is in its final position and $P_2$ has moved over to the point $P_2''$. 
Steps $2$ and $3$ involve rotations to be carried out with $P_1'$ fixed. These operations are shown in Fig. \ref{fig03}. We represent these rotations by $R_{AB}$ and
$R_{P_1'P_2'}$. $R_{AB}$ is rotation by an angle $\phi$ about the axis $AB$ (the axis perpendicular to vectors $\overrightarrow{P_1' P_2''}$ and $\overrightarrow{P_1' P_2'}$ and passing 
through the point $P_1'$) that is involved in step-$2$ above. The value of $\phi$ can vary between $0$ and $2\pi$. $R_{P_1'P_2'}$ is rotation by an angle $\theta$ 
about an axis connecting points $P_1'$ and $P_2'$ which corresponds to step-$3$ above. $\theta$ can take values from $-\pi$ to $\pi$. The sphere shown in the figure has 
got $P_1'$ at its center and has a radius equal to the distance between points $P_1'$ and $P_2''$. For convenience we have oriented the figure such that the $AB$ 
axis is vertical.

\begin{figure}
\begin{center}
\includegraphics[width=0.7\linewidth]{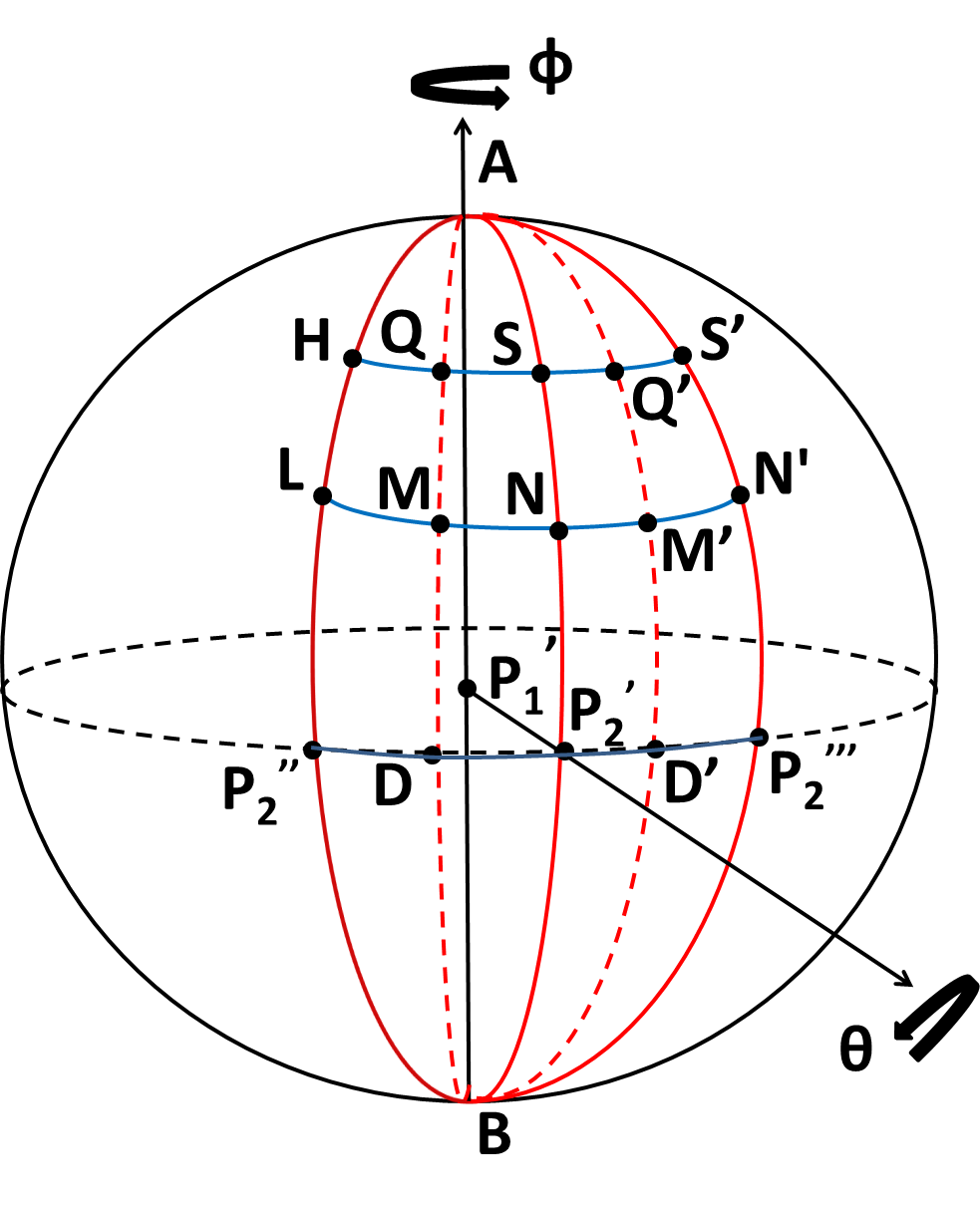} 
\caption{The figure shows how arcs of latitudes shift under the rotation about $AB$ axis. The equatorial arc $P_2''DP_2'$, $D$ being the midpoint of the arc,
moves over to $P_2'D'P_2'''$. The arc $LMN$, $M$ being the midpoint of the arc, moves over to $NM'N'$. Similarly, point $Q$, which is the midpoint of 
arc of latitude $HQS$ goes to point $Q'$. Points like $D'$, $M'$ and $Q'$ can move back to their original positions under the rotation about $P_1'P_2'$ (see Fig \ref{fig05}).}
\label{fig04}
\end{center}
\end{figure}

As the rotation $R_{AB}$ is carried out, the great circle arc $AP_2''B$ will move over into the great circle arc $AP_2'B$. And the entire region that lies between
these two arcs before rotation will now lie between the great circle arcs $AP_2'B$ and $AP_2'''B$ (see Fig. \ref{fig03}). In particular, the great circle arc $ADB$ which 
bisects the region $AP_2''BP_2'A$ (see Fig. \ref{fig04}) will move over to great circle arc $AD'B$. Consider now an arc of latitude like $LMN$, where 
$L$ lies on $AP_2''B$, $N$ on $AP_2'B$ and $M$ on $ADB$. Note that $M$
is the midpoint of the arc. Under rotation $R_{AB}$, $LMN$ will move over to latitude arc $NM'N'$. Similarly the latitude arc $HQS$ ($Q$ being the midpoint)
will move over to $SQ'S'$ under rotation (Fig. \ref{fig04}). It is interesting to note what happens to points like $M'$ and $Q'$ under the second rotation (Step-$3$ above). 
They are candidates for points that could fall back to their original position after the two rotations! 
This is so because the great circle arc $P_2'M$ ($P_2'Q$) is equal in magnitude to the great circle arc $P_2'M'$ ($P_2'Q'$) and hence under rotation about an 
axis passing through $P_2'$ and $P_1'$ both the points will fall on the same latitude circle with $P_2'$ as the pole. 

\begin{figure}
\begin{center}
\includegraphics[width=0.7\linewidth]{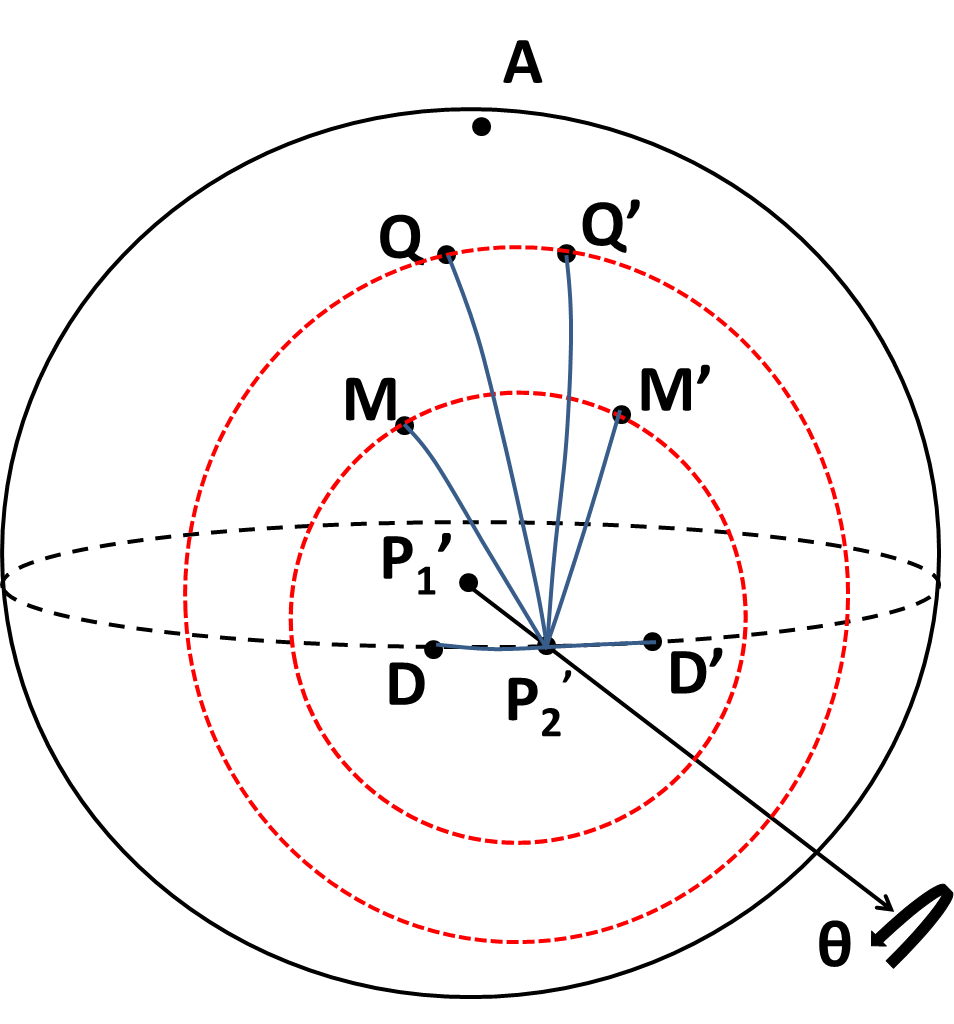} 
\caption{The pairs like points $D'$ and $D$,  $M$ and $M'$ and $Q$ and $Q'$ are equidistant from point $P_2'$. This implies that these pairs of points will
lie in the same latitudinal circles (red dashed curves) with $P_2'$ as the pole. This in turn makes it possible for these points to move back to their original 
position (that is, the one prior to rotation about $AB$) after the rotation about $P_1'P_2'$. In fact, one can show that (see text) exactly one of these set 
of points will fall back on to the original position for $\theta$ lying between $0$ and $\pi$. Note that the blue curves in the figure are great circle arcs connecting
the points involved.}
\label{fig05}
\end{center}
\end{figure}

We will now argue that depending on the value of $\theta$, there is going to be 
exactly one such point that will go back to its original position (that is, the position before Step-$2$). Fig. \ref{fig05} shows a few of the candidate points that can
come back to their original location. It is clear from the figure that the angles shown have the following ordering:
\begin{equation}
DP_2'D' = \pi > .. MP_2'M' ..>.. QP_2'Q' .. > .. AP_2'A = 0 \nonumber
\end{equation}
Even though this ordering is apparent from the figure, one can establish it more rigorously in the following manner. The spherical triangle $P_2'DQ$ (see Fig. \ref{fig04})
has in it the spherical
triangle $P_2'DM$ included. This is because the base $P_2'D$ is common for both the triangles and the great circle arcs $DQ$ and $DM$ are part of the same great circle
with $DQ$ being longer than $DM$. This implies that that the angle $NP_2'M$ is larger than $SP_2'Q$ (both being defined as angles between the corresponding
great circle arcs). But $MP_2'M' = 2MP_2'N$ and $QP_2'Q' = 2 QP_2'S$. The relationship in above equation follows.

Thus for any positive value of $\theta$ in the interval from $0$ to $\pi$, one and only one of the points of the kind $M'$ that lie in the hemisphere containing point $A$ will
come back to its original position. There would be a similar point in the diametrically opposite side of the sphere. If $\theta$ were negative and lies between $0$ and $-\pi$, 
there would be a point in the lower region below the equatorial plane that would go back to its original position and a corresponding point in the diametrically opposite side. 
Thus for any given value of $\phi$ and $\theta$ there are diametrically opposite pair of points that do not change their position under steps $2$ and $3$. This implies that
the effect of both the rotations considered above should be the same as that due to a single rotation about an axis that passes through these invariant points 
and the center of the sphere ($P_1'$). Since the effect of any arbitrary set of rotations with $P_1'$ fixed can be described using steps $2$ and $3$ above, we see that
the net effect of these rotations can be attained by a single rotation about an axis. This proves Euler rotation theorem. If we know the value of $\theta$, we can
find the invariant point by construction. To find this, pick the point (say $X$) on the great circle arc $ADB$ such that the angle between the great circle arcs
$P_2'X$ and $P_2'A$ is $\frac{\theta}{2}$. 

\section {Proof of Chasles' Theorem}
\label{sec4}
\begin{figure}
\begin{center}
\includegraphics[width=0.7\linewidth]{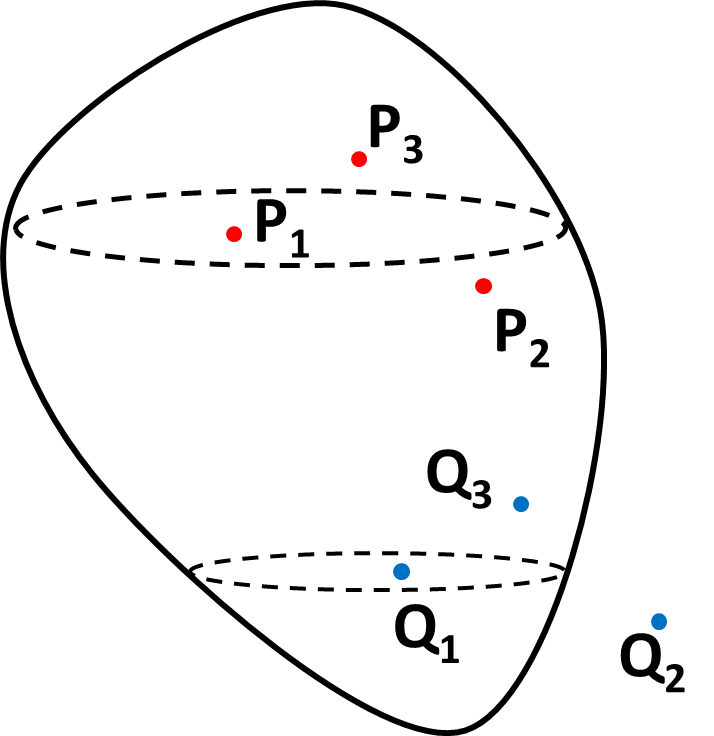} 
\caption{The proof of Chasles' theorem also involves proving that irrespective of the translating point the direction of the rotation axis and the angle of rotation
are the same. If one chooses the translating point to be $Q_1$ instead of $P_1$, the constancy of direction of axis of rotation and the amount of rotation
can be seen by considering points $Q_2$ and $Q_3$ that are related to $Q_1$ as $P_2$ and $P_3$ are to $P_1$.}
\label{fig06}
\end{center}
\end{figure}

We have already shown that the last two steps in our procedure to represent a general displacement of a rigid body corresponds to a rotation about a single axis. 
But step $1$ involved a pure translation that took point $P_1$ to $P_1'$. Thus we can conclude that a general displacement of the rigid body can be obtained by a 
translation and a rotation about an axis. To complete the proof of Chasles' theorem we also need to show that a different choice of point (instead of $P_1$) for translation will not 
alter the direction and amount of rotation to be carried out.To prove this, imagine we had chosen a different point $Q_1$ instead of $P_1$ for translation. 
Pick points $Q_2$ and $Q_3$ such that $\overrightarrow{P_1P_2} = \overrightarrow{Q_1Q_2}$ and $\overrightarrow{P_1P_3} = \overrightarrow{Q_1Q_3}$,
as shown in Fig. \ref{fig06}. One can now repeat the arguments above for proving the Euler rotation theorem. Since the vectors involved in steps $2$ and $3$ 
($\overrightarrow{Q_1Q_2}$ and $\overrightarrow{Q_1Q_3}$) in this case are identical to the old ones (even though the new displacement 
vector could be different), we will end up with the same rotation axis and angle. This completes the proof of Chasles' theorem. It may well be that there is 
no material point in the rigid body at the location of $Q_2$ or $Q_3$. One can nevertheless think of imaginary points rigidly attached to the body and moving in 
accordance with the rigidity constraints. In fact, the displacing points (like $P_1$ or $Q_1$) themselves need not form the part of the rigid body.

\begin{figure}
\begin{center}
\includegraphics[width=0.9\linewidth]{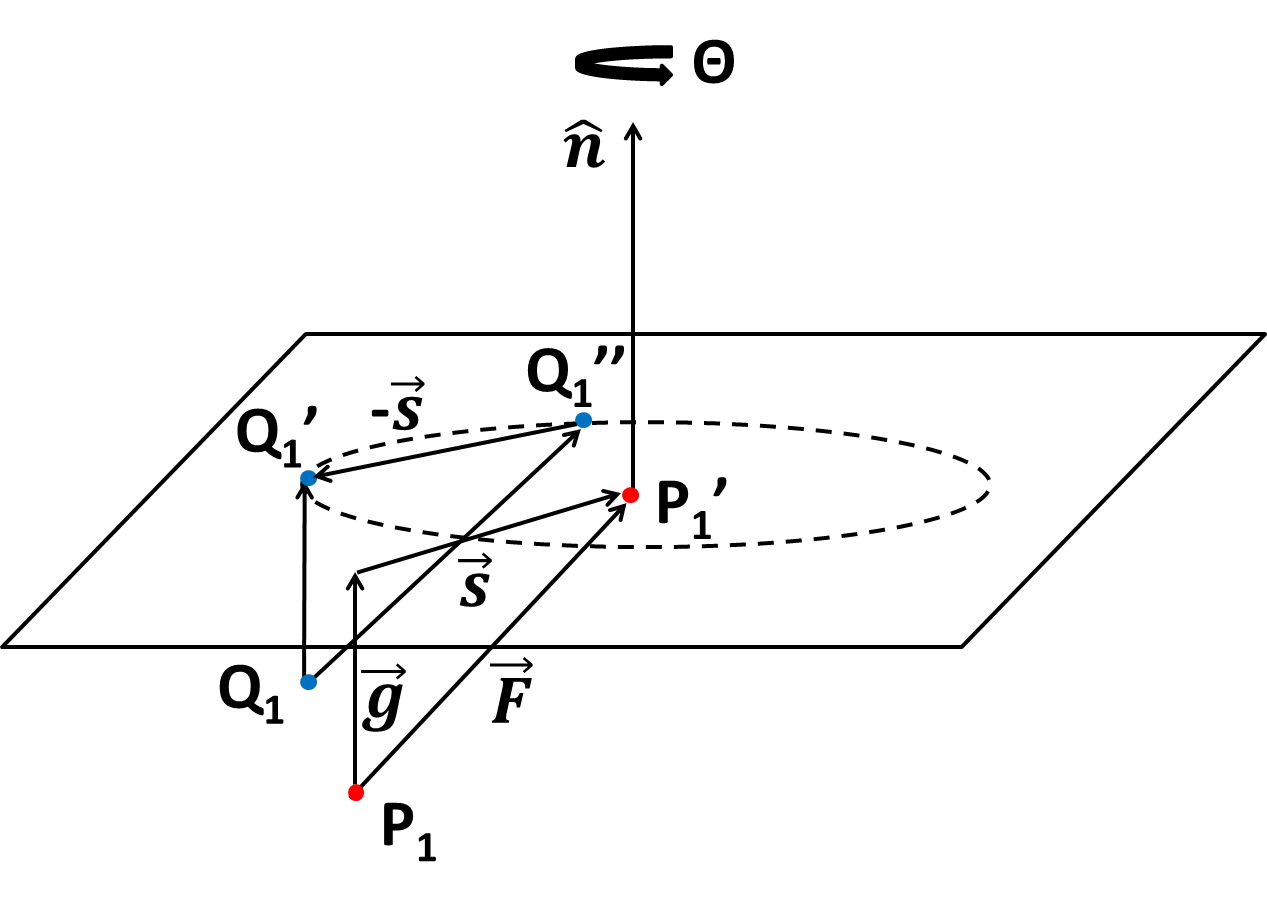} 
\caption{To prove Mozzi-Chasles' theorem, the translation vector connecting $P_1$ to $P_1'$,$\vec F$, is resolved into a part that is along the axis of rotation ($\vec g = g\hat n$) 
and another that lies the plane perpendicular to it ($\vec s$). Under an arbitrary rotation $\Theta$ about $\hat n$, one can find points that will have a 
displacement equal to $-\vec s$ (for example, displacement of $Q_1''$ to $Q_1'$ in the figure). Choosing one of these points as the translating point, one can 
arrive at a translation that is along $\hat n$. In the figure shown, picking $Q_1$ as the translating point will ensure that the translation is along the rotation axis itself.}
\label{fig07}
\end{center}
\end{figure}

An important corollary of Chasles' theorem is Mozzi-Chasles' theorem which states that a general displacement of the rigid body can be obtained by a rotation about an axis
and a translation along the same axis. For completeness, we give here a proof of this theorem. Consider a rigid body displacement described by the displacement 
vector $\vec F$ and the rotation about direction $\hat n$ by an amount $\Theta$ as shown in Fig. \ref{fig07}. $\vec F$ may not be parallel to $\hat n$.
We will  assume that the rigid displacement affects all the points in space and not just those belonging to the rigid body. 
The translation vector $\vec F$ can be expressed as a sum of vector pointing along $\hat n$ 
($\vec g$ in the figure) and a vector lying in the plane perpendicular to $\hat n$($\vec s$ in the figure). Under the rotation about $\hat n$, the different points in space 
will undergo displacements that lie in a plane perpendicular to $\hat n$. For any given value of $\Theta$, the set of displacement vectors will contain all possible vectors in the plane.
This is because the magnitude of the rotation vector will vary from zero to infinity as the distance of the points from the axis of rotation increases from zero to infinity and all points
lying on a circle at fixed distance from the axis of rotation will generate displacements in all possible directions in the plane. 
Thus one should be able to find points whose displacement is $-\vec s$ under the rotation.
If one choose one of these points as the translating point, it will ensure that the displacement vector is along $\hat n$ itself proving the corollary. The common axis
that is involved in this description, in the direction of $\hat n$, is referred to as the screw axis or Mozzi axis.

Mozzi-Chasles' theorem leads to another important result concerning motion of a rigid body in two dimensions. The counterpart of Chasles' theorem in two dimensions, sometimes
referred to as the first Euler rotation theorem, states that any displacement of a rigid body in two dimensions can be achieved by either a single rotation or a translation. 
There exists a straight forward geometric proof by construction for this theorem \cite{banach}. We shall prove the result using Mozzi-Chasles' theorem. For this, 
note that in two dimensions the axis of rotation is always perpendicular to the plane. By Mozzi-Chasles' theorem (since two-dimensional displacements are a subset of 
possible displacements in three dimensions), the rigid displacement can be achieved using a translation along an axis and rotation about that axis. 
Since the only possible translation along rotation axis is one with zero magnitude, there must be a point that does not change its position under rigid displacement 
in two dimensions. The other possibility is a pure translation in the plane in which case the screw axis will lie in the plane and the rotation about the screw axis will be zero. 
It follows that in two dimensions a rigid displacement is either a pure translation (screw axis lies in the plane) or a pure rotation (screw axis is normal to the plane).

\section{Conclusion}
\label{sec5}
We have derived Euler's rotation theorem using a novel geometric proof. The proof involves using a set of three steps that takes the
rigid body from its initial to final state. The Euler rotation theorem is derived using the last of the two steps in this procedure. The proof is presented in a 
manner that helps one in the visualization of how the invariant points arise and will be of pedagogic interest. But it should be kept in mind that the 
sequences by which one chooses to move from one configuration to the other is neither unique nor special.  
The first part of the Chasles' theorem, which asserts that the general displacement of a rigid body is a combination of translation and a rotation about 
an axis, follows immediately from Euler theorem and the first step in the procedure for carrying out rigid displacement. The fact that the orientation of axis of rotation 
and amount of rotation is independent of the translation vector involved is proved by a separate construction. For completeness, we have also presented proofs for the existence of
screw axis for motion in three dimensions and that in two dimensions any rigid displacement can be achieved by a pure rotation or a translation.
 
\begin{acknowledgements}
The author thanks Vibhu Mishra for useful discussions on existing proofs of Euler's rotation theorem.
The author would like to acknowledge financial support under the DST-FIST scheme (SR/FST/PS-1/2017/21).
This is a pre-print of an article published in Mathematical Intelligencer. The final authenticated version is available online 
at: https://doi.org/10.1007/s00283-020-09991-z
\end{acknowledgements}


\begin{thebibliography}{}
\bibitem{mozzi} G. Mozzi, Discorso matematico soprail rotamento momentaneo dei corpi. Napoli:Stamperia di Donato Campo, (1763)
\bibitem{chasles} M. Chasles, Note sur les propri\'et\'es g\'en\'erales du syst\'eme de deux corps semblables ent\'reux,  Bulletin des Sciences, Math\'ematiques,
Astronomiques, Physiques et Chemiques, 4, 321--326 (1830)
\bibitem{jackson} D. Jackson, The Instantaneous Motion of a Rigid Body, The American Mathematical Monthly, 49(10), 661-667 (1942)
\bibitem{whittaker} E. T. Whittaker, A Treatise on analytical dynamics of particles and rigid bodies. Cambridge University Press, (1947)
\bibitem{pars} L. A. Pars, A treatise on analytical dynamics. John Wiley \& Sons, (1965)
\bibitem{banach} S. Banach, Mechanics. Warszawa - Wroclaw (1951)
\bibitem{ball} R. S. Ball, The theory of screws: A study in the dynamics of a rigid body. Hodger, Foster and Company, (1876)
\bibitem{davidson} J. K. Davidson and K. H. Hunt, Robotics and screw theory: Application of kinematics and statics to robotics. Oxford University Press, (2004)
\bibitem{murray} R. M. Murray, Z. Li and S. S. Sastry, Mathematical introduction to robotic manipulation. CRC Press (1994)
\bibitem{euler} L. Euler, Formulae generales pro translatione quacunque corporum rigidorum, Novi Comm Acad. Sci. Petropol, 20, 189--207 (1776)
\bibitem{palais1} B. Palais, R. Palais and S. Rodi, A disorienting look at Euler's Theorem on the axis of a rotation, The American Mathematical Monthly, 116, 892--909 (2009)
\bibitem{palais2} B. Palais and R. Palais, Euler's rotation theorem: The axis of a rotation, J. Fixed Point Theory Appl., 2, 215--220 (2007)
\bibitem{beatty} M. F. Beatty, Kinematics of finite, rigid-body displacements, Am. J. Phys. ,34, 949--954 (1966)
\bibitem{goldstein} H. Goldstein, C. P. Poole and J. L. Safko, Classical Mechanics (3rd ed.) 134--135. Addison-Wesly, (2001)
\bibitem{thurnauer} P. G. Thurnauer, Kinematics of finite rigid-body displacements, Am. J. Phys., 35, 1145--54 (1967)
\end{thebibliography}


\end{document}